\documentstyle[11pt]{article}

\newtheorem{theorem}{{\sc Theorem}}
\newtheorem{proposition}{{\sc Proposition}}
\newtheorem{lemma}{{\sc Lemma}}
\newtheorem{corollary}{{\sc Corollary}}

\newtheorem{definition}{{\sc Definition}}

\newenvironment{proof}{\bigskip \noindent
         {\bf Proof.}}{ \hfill \fbox{\mbox{}} \medskip}

\newcommand{\cC}{{\cal C}}

\newcommand{\cL}{{\cal L}}

\newcommand{\cV}{{\cal V}}

\newcommand{\bR}{{\bf R}}
\newcommand{\bZ}{{\bf Z}}

\newcommand{\Zd}{{\bf Z}^{\, \! d}}

\newcommand{\bP}{{\bf P}}
\newcommand{\bE}{{\bf E}}

\newcommand{\bQ}{{\bf Q}}

\newcommand{\sip}{\beta} 

\newcommand{\philh}{J_n^{\cL_{1/2}}}

\newcommand{\philone}{J_n^{\cL_{1}}}  

\begin{document}

\title{\bf Formula for the Mean Square Displacement Exponent 
of the Self-Avoiding Walk in 3, 4 and All Dimensions}
\author{\sc Irene Hueter\footnotemark[1]} 
\date{} 

\maketitle

\renewcommand{\thefootnote}{\fnsymbol{footnote}} 
\footnotetext[1]{Department of Mathematics, University
                 of Florida, PO Box 118105, Gainesville, FL 32611-8105,
                 USA, email: {\sf hueter@math.ufl.edu.}} 
\footnotetext[2]{Date: August 7, 2001.}  
\footnotetext[3]{ 
              {\em Secondary:} 60D05, 60G57.}


\renewcommand{\thefootnote}{arabic{footnote}}

\vspace*{-0.8cm} 
\begin{abstract}
This paper proves the formula 
\[  \nu (d) = \left \{ \begin{array}{ll}
          1 & \mbox{for $ d=1$}, \\ 
          \max(\frac{1}{2}, \frac{1}{4} + \frac{1}{d})  
          & \mbox{for $ d \geq 2$}
        \end{array} \right . \]
for the root mean square displacement exponent $\nu(d)$
of the {\em self-avoiding walk} (SAW) in ${\bf Z}^d,$
and thus, resolves some major long-standing open conjectures rooted 
in chemical physics ({\sc Flory} (1949) \cite{flor}).
The values $\nu(2) =3/4$ and $\nu(4) = 1/2$ coincide with those that
were believed on the basis of heuristic and ``numerical evidence".
Perhaps surprisingly, there was no precise conjecture 
in dimension $3.$ Yet as early as in the 1980ies, 
Monte Carlo simulations produced a couple of confidence 
intervals for the exponent $\nu(3).$  
This work is a follow-up to {\sc Hueter} \cite{hue1}, which  
proves the result for $d=2$ and lays out the fundamental
building blocks for the analysis in all dimensions.
We consider (a) the point process of self-intersections defined via 
certain paths of length $n$ of the symmetric simple random walk 
in ${\bf Z}^d$ and 
(b) a ``weakly self-avoiding cone process'' relative to this
point process in a certain ``shape". 
The asymptotic expected distance of the
process in (b) can be calculated rather precisely as n
tends large and, if the point process has {\em circular} shape,
can be shown to asymptotically equal (up to error terms) the 
one of the {\em weakly} SAW with parameter $\sip >0.$ 
From these results, a number of distance 
exponents are immediately collectable for the SAW as well. 
Our approach invokes the {\em Palm distribution} of the 
point process of self-intersections in a cone.
\end{abstract}

\pagestyle{myheadings}
\thispagestyle{plain}
\markboth{\sc Irene Hueter}
{\sc Displacement Exponent of the Self-Avoiding Walk in All Dimensions}
                                       

\section{Introduction}

This paper is a follow-up to {\sc Hueter} \cite{hue1} and 
establishes a formula for the root mean square displacement
exponent of the self-avoiding walk in the $d$-dimensional hypercubic
lattice for all $d \geq 1.$
This simple formula resolves the ``puzzling" case $d=3$ and confirms 
the long-standing open conjectures for $d=2$ and $d=4$ that originate
in the work of {\sc Flory} \cite{flor} in the 1940ies. 
The self-avoiding walk serves as a model for linear polymer molecules.
Polymers are of interest to chemists and physicists and are
the fundamental building schemes in biological systems. A polymer
is a long chain of monomers which are joined to one another
by chemical bonds. These polymer molecules arrange themselves randomly
with the restriction of no overlap. This repelling force drives
the polymers to be more diffusive than a simple random walk.
Numerous stones wait to be uncovered from a mathematically
rigorous point of view since very little is known about the 
$2$-, $3$- and $4$-dimensional polymers or self-avoiding walk. 
At the other end, though, there is an unnumbered set of simulations and
heuristic arguments and a zoo of ``numerical artifacts" that 
lend themselves to a landscape of conjectures.
The literature has devoted much attention to this theme.
We refer the reader to other references for an overview
(e.g.\ consult {\sc Madras and Slade} \cite{masl}).
 
This paper presents some answers to the question on the  
average distance between the two ends 
of a long polymer. Our results, pertaining to the asymptotic 
expected distance of the weakly self-avoiding walk 
from its starting point up through a large step size, 
cover all \mbox{dimensions $d.$}

\smallskip
{\bf (Weakly) Self-Avoiding Walk.}
Consider the weakly self-avoiding walk in ${\bf Z}^d$ starting 
at the origin. More precisely, if $J_n= J_n(\cdot)$ denotes the
{\em number} of {\em self-intersections} or the 
{\em self-intersection local time} (SILT) 
of a symmetric simple random walk 
$S_0={\bf 0}, S_1, \ldots, S_n$ in the $d$-dimensional lattice 
starting at the origin, that is,
\begin{equation}
  \label{selfintersections}
      J_n = J_n(S_0, S_1, \ldots, S_n) = \sum_{0 \leq i < j \leq n} 
            1_{\{ S_i = S_j \}} ,
\end{equation} 
and if $\sip \geq 0$ denotes the {\em self-intersection parameter,}
then the {\em weakly self-avoiding walk} is the stochastic 
process, induced by the probability measure 
\begin{equation}
  \label{intersectmeasure}
    \bQ^{\sip}_n(\cdot) = \frac{\exp \{ - \sip J_n(\cdot) \}}{
                              \bE \exp \{ - \sip J_n(\cdot) \}} \, ,
\end{equation}
where $\bE$ stands for the expectation relative to the random walk.
In other words, $J_n = r$ self-intersections are penalized by the
factor $\exp \{ - \sip r \}.$ The measure $ \bQ^{\sip}_n$ 
may be looked at as a measure on the set of all 
simple random walks of length $n$ which weighs relative to
the number of self-intersections. This restraint walk is 
also being called the {\em Domb-Joyce model} in the 
literature (see {\sc Lawler} \cite{lawl}, p.\ 170)
but differs from the discrete {\em Edwards model}, 
which is a related repelling walk (see {\sc Madras
and Slade} \cite{masl}, 
p.\ 367 and {\sc Lawler} \cite{lawl}, p.\ 172 for some background). 
While when setting $\sip=0$ we recover the simple random 
walk (SRW), letting $\sip \rightarrow \infty$ well mimics
the {\em self-avoiding walk} (SAW). The SAW 
in ${\bf Z}^d$ is a SRW-path of length $n$
{\em without} self-intersections. 
Thus, this walk visits each site of its path {\em exactly} once.

We shall investigate the expected distance of the weakly
SAW from its starting point after $n$ steps,
as measured by the Euclidean length
and the root mean square displacement
at the $n$-th step. Let $\bE_{\sip} = \bE_{\bQ^{\sip}_n}$
denote expectation under the measure $\bQ^{\sip}_n,$ 
that is, expectation wrt.\ to the weakly SAW. Thus, 
$\bE_0$ denotes expectation wrt.\ to the SRW. Also, write  
$S_n = (X^1_n, X^2_n, \ldots, X^d_n)$ for every integer $n \geq 0.$ 
Objects of interest to us are the expectation $\bE_{\beta}$
of the {\em distance} 
$$
  \chi_n =  \{ \, \sum_{ k=1}^d (X^k_n)^2 \, \}^{1/2} 
$$ 
of the walk from the starting point ${\bf 0},$  
the {\em mean square displacement} 
$\bE_{\beta} \chi_n^2,$ and the 
{\em root mean square displacement} 
$ ( \bE_{\beta} \chi_n^2 )^{1/2}$
of the weakly SAW. Shorter, we shall write MSD and RMSD (for the
latter two), respectively. 
The RMSD exponent of the weakly SAW and the SAW, respectively,
may be defined by 
\begin{eqnarray}
    \label{rmsdexp}
      \nu_{\sip} (d) & = & \lim_{n \rightarrow \infty}
       \, \frac{ \ln  \bE_{\beta} ( \chi_n^2) }{ 2 \, \ln n} \, \\
  \nu_{\infty} (d) & = & \nu(d) =  
         \lim_{n \rightarrow \infty} \, \lim_{\sip \rightarrow \infty}
       \, \frac{ \ln  \bE_{\beta} ( \chi_n^2) }{ 2 \, \ln n} 
\end{eqnarray}
if the limits exist (otherwise we may regard the upper and lower 
exponents via $\limsup$ and $\liminf$).
Moreover, define the numbers
\begin{eqnarray}
    \label{valueofexponent}
     \mu = \mu (d) & = & 1  
          \qquad \qquad \qquad \qquad \mbox{ for } d =1, \nonumber \\
            & = &   \max( \frac{1}{2}, \frac{1}{4} + \frac{1}{d})
                \qquad \mbox{ for } d \geq 2. 
\end{eqnarray}     
Written out, $\mu(\cdot)$ takes the values 
$$
   1, \; \;  3/4, \; \;  7/12, \; \;  1/2, \; \;  1/2, \ldots.
$$      
Next, we state our main results.

\begin{theorem}
 \label{sawdistance}
The exponents of the distance of the 
weakly self-avoiding walk with $\sip >0$ and of the
self-avoiding walk in ${\bf Z}^d$ for $d \geq 1$ equal $\mu(d).$
Furthermore, there are some constants $0 < \rho_1(d) = \rho_1(d, \sip)
\leq \rho_2(d) = \rho_2(d, \sip) < \infty$  such that 
$$  
   \rho_1(d) \leq
     \liminf_{n \rightarrow \infty}
       n^{- \mu(d)} \, \bE_{\sip} (\chi_n) 
        \leq  \limsup_{n \rightarrow \infty}
       n^{- \mu(d)} \, \bE_{\sip} (\chi_n) 
         \leq \rho_2(d) ,
$$
where $\rho_1(d)$ is uniform in $\sip$ for $d \geq 5$ and 
may depend on $\sip$ for $d \leq 4$ and
$\rho_2(d)$ is uniform in $\sip$ for $d \leq 2$ and $d \geq 5$ and 
may depend on $\sip$ for $d =3,4.$ 
\end{theorem}

The proof is in Corollary \ref{liminfsupmsd} for $d \geq 2$ and 
in Theorem \ref{onedimension} for $d=1.$
 
\begin{theorem}
  \label{rmsd}
The weakly self-avoiding walk with $\sip > 0$ and 
the self-avoiding walk in ${\bf Z}^d$ for $d \geq 1$ have
$$ 
    \nu_{\sip} (d) = \nu(d) = \mu(d).
$$
Moreover, there are some constants  $0 < \rho_3(d) = \rho_3(d, \sip)
\leq \rho_4(d) = \rho_4(d, \sip) < \infty$  such that 
$$
   \rho_3(d)  \leq
     \liminf_{n \rightarrow \infty}
       n^{- 2 \mu(d)}  \, \bE_{\sip} (\chi_n^2) 
        \leq  \limsup_{n \rightarrow \infty}
       n^{- 2 \mu(d)}  \, \bE_{\sip} (\chi_n^2) 
         \leq \rho_4(d),
$$ 
where  $\rho_3(d)$ is uniform in $\sip$ for $d \geq 5$ and 
may depend on $\sip$ for $d \leq 4$ and
$\rho_4(d)$ is uniform in $\sip$ for $d \leq 2$ and $d \geq 5$ and 
may depend on $\sip$ for $d =3,4.$ 
\end{theorem}

See Corollary \ref{variance} for a proof when $d \geq 2$ and
Corollary \ref{vardimone} when $d=1.$
{\sc Hueter} \cite{hue1} proves the analogous results 
in the two-dimensional context.  
Theorem \ref{rmsd} settles a couple of major decades-old 
open conjectures that can be traced back to at least
{\sc Flory}'s work \cite{flor} in the 1940ies.
It is hoped that our point of view will shed light onto the 
(weakly) SAW enough to bring to fruition solutions of
other tantalizing problems on these and related objects.
Whereas our results here and in \cite{hue1}
are novel for $d=2,3,$ and $4$ and 
$\sip \in (0, \infty],$
the result on the RMSD exponent for the SAW for $d \geq 5$ 
is in {\sc Hara and Slade} \cite{hasl0, hasl}
and the one on the RMSD exponent for the weakly SAW for $d=1$ 
in {\sc Greven and den Hollander} \cite{grho}. Of course,
the result on the one-dimensional SAW is quite obvious. 
The former was investigated via the perturbation technique
``lace expansion" and the latter via large deviation theory.
{\sc Brydges and Spencer} \cite{brsp} establish that the
scaling limit of the weakly SAW is Gaussian
for sufficiently small $\sip >0$ and $d \geq 5.$

A couple of Monte Carlo simulations were performed as early as
in the 1980ies to estimate
the RMSD exponents for the SAW (for more references and details
on this, see {\sc Madras and Sokal} \cite{maso}). The produced 
$95 \%$-confidence intervals appear to center around the 
value $0.59....$ and would suggest a value slightly larger
than $7/12=0.58333...$ On another historical note,
an earlier estimate was the Flory estimate $0.6.$
Just for $d=4,$ a {\em logarithmic correction} associated with 
$\bE_{\beta} \chi_n^2$ is being predicted
(visit e.g.\ {\sc Lawler} \cite{lawl}, p.\ 167). 
We point out that our results leave space for such a correction
for $d=3, 4.$ Indeed, the expressions that we derive for both
the mean square displacement and the expected distance of
the weakly SAW in this paper are bounded by constants 
that depend on $\sip$ as $\sip \rightarrow \infty.$ 
In order to exchange limits as $\sip \rightarrow \infty$
and as $n \rightarrow \infty$ of $\bE_{\beta} (\chi_n^2) /
n^{2 \mu(d)}$ as is necessary to obtain the MSD of the
SAW, we would need to know how $\sip$ and $n$ are related,
in other words, in how far $\sip$ is bounded by some function in
$n$ and vice versa. Such a relationship would allow to 
translate the bounds of $\bE_{\beta} (\chi_n^2) /
n^{2 \mu(d)}$ in terms of $\sip$ into bounds in $n,$ 
and thus, into some correction factors to $ n^{2 \mu(d)}.$ 
Nevertheless, none of this is needed to extend the
values of the distance and MSD exponents of the weakly
SAW to the SAW. Also, while our results are not new for $d=1$ and
$d \geq 5,$ our approach provides an alternative proof.

Perhaps surprisingly, {\em one} and the same approach -- the
one employed in this work -- suffices to handle all dimensions 
$d \geq 1.$ With a bit of extra work, the results for $d \geq 3$
follow from the analysis for the case $d =2,$ whereas the case $d=1$
has a different touch. Therefore, the latter dimension is dealt 
with in a separate section (Section 6).
In dimension $4,$ an interesting twist occurs to the expression
for the expected distance of the (weakly) SAW, and thus, for
the distance exponent as well. This expression carries two
significant terms, one of which is dominating in dimensions
$2$ through $3,$ the other of which is dominating in dimensions
$ d \geq 5$ and may be identified as the term that resembles
the contribution which we would obtain for the SRW. In dimension
$4,$ both terms compete with each other. While the SRW-term wins 
for all $\sip >0$ below a certain threshold,
it is not clear which of both terms dominates for large $\sip.$
In this sense, for instance, dimension $4$ is more intriguing than 
dimension $3$ despite the fact that 
the RMSD exponent is the same as for the SRW.
Hence, in dimension $4,$ the exponent $1/2$ arises
for different reasons than it occurs in dimensions $5$ and
higher. 

Our strategy of proof is to regard a process 
which is penalized according to the number of self-intersections 
of the random walk that we see in {\em one direction} and 
to compare its expected distance to the one of the weakly SAW.
For this purpose, we will spread a fixed collection of rays
that emanate from the origin, each of which describes a cone.
Some of these cones will carry a more typical number of
self-intersections than others -- typical will mean of 
order $\sqrt{n}.$ Furthermore, the event that a cone is less typical
will depend on the realized SRW-path. For $d \geq 3,$ 
the space becomes large in the sense that
the cones which contribute most of the self-intersections
of the weakly SAW have cardinality of order {\em strictly}
less than the order of the total number of 
cones (see (\ref{lowsiltlines})). 

We will invoke the {\em Palm distribution}
of the point process of self-intersections, defined
via certain paths of length $n$ of the symmetric
SRW in $\Zd,$ in a cone to introduce 
a ``weakly self-avoiding cone process'' relative 
to this point process when in a certain ``shape''.
The asymptotic expected distance of this
process can be calculated rather precisely as n
tends large and, if the point process has {\em circular} shape,
can be shown to asymptotically equal (up to error terms) the 
one of the weakly SAW with parameter $\sip >0.$
Then we collect analogous results on the 
mean square displacement of the weakly SAW.
From these results, upon some considerations
on uniform bounds and estimates in $\sip$ as 
$\sip \rightarrow \infty,$ the distance and the MSD exponents
of the SAW immediately derive.

The paper is organized as follows.
Section 2 specifies the SRW-paths that are
significant from a weakly SAW's point of view.
Section 3 makes a connection between Palm distributions 
and the random walk and recalls the notions of shape of the
underlying point process and of
a weakly self-avoiding cone process.
Section 4 calculates some asymptotic 
mean distances of this process and links those to the ones 
of the weakly SAW. Section 5 discusses the transfer of the
distance and MSD exponents to the SAW.
Some remarks on the transitions
$\sip \rightarrow \infty$ and $\sip \rightarrow 0$ end 
Section 5. Finally, Section 6 takes care of the one-dimensional
setup.


\section{SILT that is Typical for the Weakly SAW}
\setcounter{equation}{0} 

Most of the remainder of the paper will be devoted to studying
the weakly SAW. We shall exploit the information that 
is contained in the intersections that are discouraged
but not forbidden as for the SAW. In low dimensions, 
the weakly SAW pays attention to the SRW-paths that exhibit 
a smaller number of self-intersections than is expected for the
SRW. This effect is most emphasized in dimension $1.$ Paths
that have about $\bE_0 J_n$ self-intersections are not
important from the perspective of a weakly SAW. While
a weakly SAW-path of length $n$ will turn out to have 
expected SILT of order $n$ in all dimensions, the SRW is 
forced to intersect itself more frequently, 
at least in dimensions $1$ and $2.$
We begin to review the average $\bE_0 J_n$ for the SRW 
and to derive the range for $J_n$ that is significant
from the point of view of the weakly SAW.

A favorite exercise in a probability course is as follows.
By invoking the Fubini theorem and the Local Central Limit 
theorem, we obtain for all sufficiently large even $n,$
\begin{eqnarray}
  \label{srwsilt}
 \bE_0 J_n & = & \sum_{0 \leq i < j \leq n} \bP_0 ( S_i = S_j ) 
                 \\[0.15cm] 
           & = &  (1 + o(1)) \, \sum_{0 \leq i < j \leq n/2 }
                       2 \, (\frac{d}{2 \pi (j-i)})^{d/2}
                       \nonumber 
\end{eqnarray} 
       \[  \mbox{} \hspace*{0.5cm} 
             =  (1+o(1))  \left \{  \begin{array}{ll}
                           ( \frac{2}{3 \, \pi^{1/2}} \, n^{3/2} 
                                     & d=1, \\
                            \frac{1}{\pi} \, n \ln n  
                                      & d=2, \\
                             c_d \, n          & d \geq 3 
                            \end{array} \right.  \]         
for some positive finite constants $c_d,$
where we used the $o(\cdot)$ notation, that is,
write $f(n)= o(g(n))$ as $n \rightarrow \infty$ for two real-valued
functions $f$ and $g$ if $\lim_{n \rightarrow \infty} f(n)/g(n) =0.$

Now, let $\omega_0 = \omega_0(d)$ denote the logarithm of the
{\em connective constant} or the exponent of the
number of SAW-paths, in other words, the exponential rate 
at which the cardinality of the set of all SAW-paths 
$S_0={\bf 0}, S_1, \ldots, S_n$ (with $J_n =0$) up through time $n$
grows in $n.$ Observe that $\omega_0(1) = 0$ and
$ d \leq  e^{\omega_0(d)} \leq 2d-1$ for all $d.$
The {\em upper} bound $2d-1$ for $e^{\omega_0}$ may be
seen by counting all paths of length $n$ that do not return
to the  most recently visited point (clearly, an overestimate),
whereas the {\em lower} bound $d$ for $e^{\omega_0}$ may be seen
by counting all paths of length $n$ that take only positive
steps in both coordinates, for example for $d=2,$
i.e.\ move only north or east, say.

The two subsequent Propositions restate the results in 
Propositions 1 and 2 in 
{\sc Hueter} \cite{hue1}, Section 2, without proofs.
Since the arguments of proof are carried out in the time
space, as opposed to the state space, they as well apply for
$d \not = 2.$ 
When reading the proofs of {\sc Hueter} \cite{hue1}, written
for $d=2,$ the reader might want to replace the number `$4$',
the number of nearest neighboring sites of each lattice site, 
by $2d$ and rely on $\omega_0(d),$ as just described, rather
than $\omega_0(2).$  

The idea to prove the {\bf upper} bound is that 
it suffices to find a subset of SRW-paths that contributes 
strictly more to  $ \bE_0 \exp \{ - \sip J_n  \}$ than
the set of paths with $J_n > Bn$ for all $B > B_*$ and
some suitable positive finite constant $B_*.$
In fact, the set of all self-avoiding paths satisfies 
this requirement. It is enough to derive a lower
bound for 
$\bP_0 ( J_n =0 ) = \bE_0 ( \exp \{ - \sip J_n \} \, 
1_{\{ J_n =0 \}} )$ and to see for which $B_*$
it is strictly larger than 
$ \exp \{ - \sip B_* n \} >  \bE_0 ( \exp \{ - \sip J_n \} 
\, 1_{\{ J_n > B_* n \}} ).$

\begin{proposition}
   \label{upperboundforjn}{\bf (Upper Bound for $J_n$)}
Let $d \geq 1,$ $\sip > 0,$ and let $\omega_0(d)$ denote
the exponent of the number of self-avoiding walks.
Then for every $B > B_* = B_*(d) = (\ln (2d) - \omega_0(d))/\sip >0$ 
and every integer $n \geq 0,$
$$
   \bE_0 ( e^{ - \sip J_n} \, 1_{\{ J_n > B n \}} ) < 
            \bE_0 ( e^{ - \sip J_n} \, 1_{\{ J_n =0 \}} ),  
$$
in particular, as $n \rightarrow \infty,$
$$
   \bE_0 ( e^{ - \sip J_n} \, 1_{\{ J_n > B n \}} ) = 
            o( \bE_0 ( e^{ - \sip J_n} \, 1_{\{ J_n =0 \}} )) .  
$$
\end{proposition}

\begin{proof}
The proof is in {\sc Hueter} \cite{hue1}, Proposition 1.
\end{proof}

Hence, we may restrict our attention to the SRW-paths 
that exhibit $J_n \leq n B_*.$
On the other hand, the paths with $J_n $ of order
less than $n$ are not significant either.

\begin{proposition}{\bf (Lower Bound for $J_n$)} 
    \label{lowerboundforjn} 
Let $d \geq 1$ and $\sip > 0.$ There is some $\zeta_*(\sip) >0 $ 
that can be made precise
(\cite{hue1}, Proposition $2$)
such that for $b_* = \zeta_*(\sip)/\sip >0,$ 
for every $\delta >0$ and every $b < b_*,$
as $n \rightarrow \infty,$ 
$$
    \bE_0 ( e^{ - \sip J_n} \, 1_{\{ J_n \leq n^{1- \delta} \}} ) =
      o(   \bE_0 ( e^{ - \sip J_n} \, 1_{\{ J_n < b n \}} )).
$$       
\end{proposition}

\begin{proof}
The proof is carried out in {\sc Hueter} \cite{hue1}, Proposition 2.
\end{proof}

The reasoning to prove the {\bf lower} bound is that
it is sufficient to
identify a subset of SRW-paths that contributes strictly more
to   $ \bE_0 \exp \{ - \sip J_n  \}$ than the set of paths
with $J_n  \leq n^{1- \delta} .$ 
The latter set is modified
by introducing of order $n$ repetitions of steps to each
SRW-path, which gives rise to a set of paths that have
$J_n \leq bn$ for $b < b_*,$ whose size is 
of strictly larger exponential order.
Then $b_*$ can be chosen suitably small.

We remark that the same proof with slight adjustments applies
when the bound $n^{1 - \delta}$
in the statement of Proposition \ref{lowerboundforjn} 
is replaced by $n \, q_n,$ where $q_n \rightarrow 0$ arbitrarily
slowly as $n \rightarrow \infty.$  Hence, the set of paths
with $J_n \in [0, n q_n)$ contributes to $ \bE_0 ( e^{ - \sip J_n})$
or to the $k$-th moments $\bE_{\sip}(\chi_n^k)$ only in a negligible
fashion in the sense that the contribution is $o(\bE_0 ( e^{ - \sip J_n}))$
or $o(\bE_{\sip}(\chi_n^k)),$ respectively, as $n$ tends large
(in fact, this error term is exponentially smaller,
as the proof of Proposition \ref{lowerboundforjn} indicates). 
As a consequence of Propositions \ref{upperboundforjn} and
\ref{lowerboundforjn}, for all that follows, we may neglect 
to keep track of those error terms and assume that 
\begin{equation}
  \label{rangejn}
   J_n \in [ b_1 n,  b_2 n ]
\end{equation}
for all sufficiently large $n$ and 
for some constants $0 < b_1 < b_2 < \infty$ 
such that $\sip b_2$ is a positive number independent of $\sip$ 
and $\sip b_1$ may depend on $\sip$ in such a way that
$\sip b_1 $ tends to zero as $\sip \rightarrow \infty.$
Observe that comparing (\ref{srwsilt}) and (\ref{rangejn}) 
along with the observations in the last paragraph reveals that,
for $d \leq 2,$ the expectation $\bE_0 J_n$ is of larger order
in $n$ than $\bE_{\sip} J_n,$ whereas for $d \geq 3$ 
and $\sip$ below a certain threshold, $\bE_0 J_n < \bE_{\sip} J_n.$


\section{Point Process of Self-Intersections and Cones}
\setcounter{equation}{0}  

Throughout this section, we will assume that $d>1,$ and 
for the rest of the paper, we shall omit discussion of the 
obvious case $\sip =0.$ 
If we let $X^1_n,$ $ X^2_n, \ldots X^d_n$ denote the 
coordinate processes of the SRW, that is, 
$S_n = (X^1_n, \ldots , X^d_n)$ for every 
integer $n \geq 0,$ define the distance
\begin{equation}
 \label{distance}
   \chi_n =  \vert \! \vert S_n \vert \! \vert
             = \{ \, \sum_{ k=1}^d (X^k_n)^2 \, \}^{1/2}  
\end{equation}  
or the root of the square displacement
of the walk from the starting point ${\bf 0}.$  
Furthermore, let $ \bP_{\chi_n}$ denote the probability distribution
of the distance $\chi_n$ of the SRW. 
 
As in {\sc Hueter} \cite{hue1}, we shall rely on a process
which is intimately related to the weakly SAW.
Let us think about asymptotically calculating the expected 
distance of the SRW after $n$ steps
from the starting point (for which process, though, the 
calculation is more straightforward). 
One possible route involves approximating the SRW by
means of a Brownian motion and controlling the entailed errors,
which may be sketched in the following way.
Rely on the Local Central Limit theorem and 
rewrite the density of the approximating Brownian motion to the
SRW in polar coordinates. Then the asymptotic expected distance
of the SRW is calculated via integrations over the {\em radial part}
and the {\em angle.}
In case of the weakly self-avoiding process, 
the penalizing weight takes into consideration the number of
self-intersections in a direction, that is,
near the line that passes through the starting 
point and the endpoint of the SRW-path. 
Part of our strategy will consist in relating the expected distance
of this newly-defined process with the one of the weakly SAW
and in finding bounds on the expected distance of the former 
process by \\*[0.05cm]
\mbox{} \,
(a) keeping track of the radial part of the SRW, penalized by
the SILT in a certain cone,
\\*[0.05cm] \mbox{} \,
(b) by integrating out over all lines in $\cV.$  
 
\bigskip
{\bf 3.1. Point Process of Self-Intersections and Cones.}
The next subsection will utilize Palm distributions of the point
process of self-intersection points of the SRW with 
$J_n \in [ b_1 n,  b_2 n ]$ to define typical 
penalizing weights within certain classes of cones.
Palm distributions help answer questions dealing with properties
of a point process, viewed from a {\em typical} random geometric
object that is defined via the point process. 
As a simple example we could explore the mean number 
of points of a point process in the plane whose nearest 
neighbors are all at distance at least $r.$

Let us, however, first recall the notation set in 
{\sc Hueter} \cite{hue1} to describe the point process 
of self-intersections and its associated cones.
Let $\Phi= \Phi_n = \{ x_1, x_2, \ldots \}$ denote the point process
of {\em self-intersection points} of the SRW in ${\bf Z}^d$ when
$J_n \in [ b_1 n,  b_2 n ].$ Note that $ \vert \Phi \vert
\in [ b_1 n,  b_2 n ]$ and $\Phi$ depends on $n,$
$b_1,$ and $b_2,$ thus, on $\sip.$ We allow the points $x_i$ of $\Phi$
to have multiplicity and count such a point exactly as many
times as there are self-intersections of the SRW at $x_i.$
This random sequence of points $\Phi$ in ${\bf Z}^d$ 
may also be interpreted as a random measure. 
Note that $\bE_0 \Phi$ is $\sigma$-finite.
Let $N_{\Phi}$ denote the set of all point sequences, generated
by $\Phi,$ ${\cal N}_{\Phi}$ the point process $\sigma$-algebra 
generated by $N_{\Phi},$ and $\varphi \in N_{\Phi}$ 
denote a realization of $\Phi.$ Formally, $\Phi$ is a measurable
mapping from the underlying probability space into 
$(N_{\Phi},{\cal N}_{\Phi})$ that induces a distribution on
$(N_{\Phi},{\cal N}_{\Phi}),$ the distribution ${\bP_{\Phi}}$ 
of the point process $\Phi.$ By virtue of the $\sigma$-finiteness 
of $\bE_0 \Phi,$ ${\bP_{\Phi}}$ is a probability measure. 
Also, let $\bE_{\Phi}$ denote expectation relative \mbox{to 
$\bP_{\Phi}.$} For more details on this language, e.g.\ consult
 {\sc Stoyan, Kendall, and Mecke} \cite{skm}, see Chapter 4, 
p.\ 99. 

We will first calculate expected distance measures 
for a related process, called a ``weakly self-avoiding
cone process" relative to $\Phi$ when in a certain ``shape".
This process appears to us to have the advantage over the weakly SAW
that rather precise estimates can be calculated for the expected 
distance from the origin. For this purpose, our study will center
around how to distribute the self-intersection points $\Phi$ 
to cones which are positioned at the origin in $\Zd.$

A cone may be specified by a line that the cone contains.
Thus, we introduce a {\em test set} $\cV$ of half-lines $L$ -- that
we shall call {\em lines} for short -- which emanate from the origin
and the intersection points of which with the $d$-dimensional 
unit sphere are uniformly and regularly distributed over the sphere.
It will not become clear until much later that this is a 
possible optimal way of choosing the lines for $\cV$
(see, for instance, \mbox{Definition \ref{shapeofv}}
and Lemma \ref{conditiondsatisfied}).
We will postpone determining the cardinality $\vert \cV \vert$ of $\cV$
to the proofs of Propositions \ref{msdupper} and \ref{msdlower}, 
which will be the only relevant fact about $\cV$ to retain.
Next, for any $L \in \cV,$ let the ``cone" $\cC_L$ be defined by
\begin{eqnarray}
  \label{cone}
  \cC_L &  =&  \{ x_i \in \Phi : \mbox{dist}(x_i , L) \leq 
                           \mbox{dist}(x_i , L') \mbox{ for all }
                           L \not = L' \in \cV \}
\end{eqnarray}
with the convention that if equality 
 $\mbox{dist}(x_i , L) = \mbox{dist}(x_i , L')$ holds for two lines
 $L$ and $L'$ and a certain number of points $x_i,$ then half of them
 will be assigned to $\cC_L$ and the other half to $\cC_{L'}.$ 
 Note that no point of $\Phi$ belongs to more than one $\cC_L$ 
and each point to exactly one $\cC_L.$ Thus, $\vert \cC_L \vert $ 
equals the SILT of the SRW $S_n$ in a cone at the origin that contains
the line $L.$ 
Once the lines are selected for $\cV,$ we may classify them 
according to the SILT that their cones carry.
For any constants $0 < a_1 < a_2 < \infty,$ 
for any suitably small $\delta >0,$ and for each
$ 0 \leq r \leq 1,$ define the random sets 
\begin{eqnarray}
  \label{halfline}
  \cL_{1/2} & = & \cL_{1/2}(\Phi) =  
           \{ L \in \cV: \, 2 \vert \cC_L \vert \in 
              [ a_1 n^{1/2}, a_2 n^{1/2} ] \} \\
  \cL_{1/2 \pm } & = & \cL_{1/2 \pm }(\Phi)  =  
        \{L \in \cV: \, 2 \vert \cC_L \vert
               \in [ a_1 n^{1/2 - \delta},  a_2 n^{1/2 + \delta}] \} 
               \nonumber \\
   \cL_{-} & = & \cL_{-}(\Phi)  = 
           \{ L \in \cV: \, 2 \vert \cC_L \vert 
                    \in (0,  a_1 n^{1/2 - \delta}) \} \nonumber \\
    \cL_{+}  & = & \cL_{+}(\Phi)  =  
             \{ L \in \cV: \, 2 \vert \cC_L \vert
                     \in ( a_2 n^{1/2 + \delta}, 2 b_2 n] \} \nonumber \\ 
   \cL_{r} & = & \cL_{r}(\Phi) = 
              \{ L \in \cV: \, 2 \vert \cC_L \vert
              \in [ a_1 n^{r}, a_2  n^{r} ] \}  
                          \nonumber \\  
   \cL_{\emptyset} & = & \cL_{\emptyset}(\Phi) = 
              \{ L \in \cV: \, \vert \cC_L \vert = 0 \}, 
                   \nonumber
\end{eqnarray} 
which depend on $a_1, a_2,$ and $\cV.$ We will choose
$a_1$ and $a_2$ such that $a_1 \sip$ and $a_2 \sip $ 
are positive numbers which are independent of $\sip$ and $n.$

\medskip
{\bf 3.2. Weakly Self-Avoiding Cone Process relative to
$r$-Shaped $\Phi.$}
If $h: \bR \times N_{\Phi} \rightarrow \bR_+$ 
denotes a nonnegative measurable real-valued function and 
$\cL_*(\Phi)$ denotes any subset of lines
in $\cV,$ then
since $\bE_0 \Phi$ is $\sigma$-finite, we may disintegrate
relative to the probability measure $\bP_{\Phi},$ 
\begin{equation}
 \label{disintegration}
    \bE_{\Phi}  \left ( \sum_{L \in \cL_*(\Phi)}
                              h(L, \Phi) \right )
               = \int  \sum_{L \in \cL_*(\varphi)}
                              h(L, \varphi)\,  d \bP_{\Phi} (\varphi)
\end{equation} 
(visit also {\sc Kallenberg} \cite{kall}, p.\ 83, and 
{\sc Stoyan, Kendall, and Mecke} \cite{skm}, p.\ 99). 
For a discussion of some examples of Palm distributions of
$ \bP_{\Phi},$ the reader is referred to the Appendix in 
{\sc Hueter} \cite{hue1}.  

Observe that the conditional distribution
$\bP_{\Phi \vert \chi_n }$ of the 
point process $\Phi,$ given $\chi_n,$ is a function of $\chi_n$ 
and depends on condition (\ref{rangejn}), as explained earlier,
so as to produce realizations that exhibit
$J_n \in [ b_1 n,  b_2 n ].$ 
Apply formula (\ref{disintegration}) with
\begin{equation}
   \label{choosehfct}
   h(L,\Phi) = \frac{ \exp \{ - \sip \vert \cC_L \vert \}}
                    { \vert \cL (\Phi) \vert} ,
\end{equation}
with $\bP_{\Phi \vert \chi_n }(\varphi \vert x)$ in place of 
$\bP_{\Phi} (\varphi),$ and $\cL_* = \cL \subset \cL_{r} \subset
\cV $ to define the numbers $a_x= a_x(\cL)$ by
\begin{eqnarray}  
   \label{axnumbershalfr}
       \exp \{ - \sip a_x  n^{r}/2 \} & = &
          \bE_{\Phi \vert \chi_n} ( \, \vert \cL (\Phi) \vert^{-1}  
                          \sum_{L \in \cL (\Phi)} 
                     e^{- \sip \vert \cC_L \vert  } \vert 
                         \chi_n =x )  \\*[0.1cm]
                 & = & 
             \int_{\bZ^d} \, \vert \cL (\varphi) \vert^{-1}  
                \sum_{L \in \cL (\varphi)} 
                     e^{- \sip \vert \cC_L \vert } \, 
                         d \bP_{\Phi \vert \chi_n }(\varphi \vert x)
\nonumber
\end{eqnarray} 
for $0 \leq x \leq n,$ where we set $ \sum_{L \in \cL} = 0$ 
if $\cL = \emptyset.$
For example, if we set $\cL = \cL_{1/2}$ and $r=1/2,$
then conditioned on the event $\chi_n =x ,$ the number
$a_x(\cL_{1/2}) n^{1/2}/2$ may be interpreted as 
``typical" SILT relative to the lines in $ \cL_{1/2},$ 
equivalently, $\exp \{ - \sip a_x  n^{1/2}/2 \}$ represents a
``typical'' penalizing factor with respect to $ \cL_{1/2},$
provided that $\chi_n =x.$
Taking expectation, we arrive at the
expected ``typical'' penalizing factor 
\begin{equation}  
  \label{averageliner} 
    \bE_0 ( e^{- \sip J_n^{\cL}})
            =   \bE_0 ( \exp \{ - \sip a_{\chi_n}  n^{r}/2 \}) .
\end{equation}
In the same manner, we calculate  
\begin{equation}  
  \label{distaverpenweight}
   \bE_0 ( \chi_n \, e^{- \sip J_n^{\cL}} )  
    = \bE_0 ( \chi_n \, 
            \bE_{\Phi \vert \chi_n} ( \, \vert \cL (\Phi) \vert^{-1}  
                          \sum_{L \in \cL (\Phi)} 
                     e^{- \sip \vert \cC_L \vert  } \vert 
                         \chi_n =x ) ) .
\end{equation}
The proofs of Propositions \ref{msdupper} and \ref{msdlower}
below (see also Definition \ref{coneprocess}) will shed
light on the issue of this particular choice of penalizing
weight. It will turn out that a crucial role will be played by
variants of the quotient  $ \bE_0 (  \chi_n \, e^{- \sip \philh} ) / 
 \bE_0 ( e^{- \sip \philh} ).$   

\begin{definition}[{\bf $\Phi$ or $\cV$ are $r$-shaped}]
  \label{shapeofv}
Let $\rho >0$ be suitably small.
We say that $\cL_r$ {\em contributes (to $J_n$) essentially}
if
$$ 
    \sum_{L \in \cL_r} \vert \cC_L \vert \geq  \frac{1}{2} \,
                     J_n^{1- \rho}. 
$$
In this case, we say that $\cV$ and $\Phi$ are {\em $r$-shaped} or
have {\em shape $r.$}
In particular, when $r=1/2,$ then we say that $\cV$ and $\Phi$ 
have {\em circular shape} or are {\em circular.}
The convention is that multiple shapes
are allowed, that is, $\Phi$ may simultaneously 
have shape $1/2$ and shape $3/4.$
\end{definition}
{\bf Remarks.} \\
{\bf (1)} For our purposes and later calculations, 
it is not necessary that the lines contributing essentially,
as explained in Definition \ref{shapeofv}, have exact SILT
of order $n^r$ in the sense that the real value $r$ is hit
precisely. Instead, it suffices to replace $\cL_r$ by
$\cL_{r \star} =
 \{L \in \cV: \, 2 \vert \cC_L \vert
      \in [ a_1 n^r,  a_2 n^{r + \delta}] \} $ for $\delta >0,$
and to ultimately let $\delta \rightarrow 0$ in the obtained
results (because $\delta >0$ was arbitrary). Hence, when
applying Definition \ref{shapeofv}, we will think of
$\cL_{r \star} $ rather than $\cL_r$ and refer to
\begin{equation}
  \label{approxshape} 
   \sum_{L \in \cL_s \atop \mbox{for } r \leq s \leq r + \delta } 
    \vert \cC_L \vert \geq  \frac{1}{2} \, J_n^{1- \rho} . 
\end{equation} 
With this meaning, it is obvious that, 
for sufficiently large $n,$ there must be $0 \leq r \leq 1$
such that the set $\cL_{r \star} $ contributes essentially, and thus, 
the shape of $\Phi$ and $\cV$ is well-defined.
Nevertheless, for the sake of not complicating our presentation,
we shall not write $\cL_{r \star} $ and not use the extension
in (\ref{approxshape}) but simply write $\cL_r.$ \\*[0.1cm]
{\bf (2)} 
We might as well choose $J_n \, \tau_n /2$ with $\tau_n \rightarrow
0$ arbitrarily slowly as $n \rightarrow \infty$ in place of 
$J_n^{1- \rho}/2$ in the defining inequality for the shape of $\Phi.$ 
There is nothing special about the choice above. 

\smallskip
Next, if $\cL_r$ contributes essentially then,
by (\ref{rangejn}) and (\ref{halfline}), 
\begin{eqnarray}
  \label{boundsonlr}                     
   \frac{b_1}{a_2} \, n^{1-r-\rho} & \leq & 
        \vert \cL_r \vert \leq \frac{2 b_2}{a_1} \, n^{1-r} .
\end{eqnarray}
It is apparent that the upper bound in (\ref{boundsonlr})
holds even when $\Phi$ is not $r$-shaped. Since we choose 
$a_1$ and $a_2$ such that $\sip a_1$ and $\sip a_2$ are
independent of $\sip,$ it follows that $b_2/a_1$ is
independent of $\sip.$

\begin{definition}[{Weakly self-avoiding cone process relative 
to $r$-shaped $\cV$}]
  \label{coneprocess}
Define a {\em weakly self-avoiding cone process relative 
to $\cV$ in shape $r$} by some $d$-dimensional process whose
radial part is induced by the probability measure 
\begin{equation}
  \label{distconeprocess}
    \bQ^{\sip, \cV, r}_n = \frac{\exp \{ - \sip \vert \cC_L \vert \}}{
                       \bE_0 \exp \{ - \sip  J_n^{\cL_r}  \}} 
\end{equation}
on the set of SRW-paths of length $n$
if $\cV$ has shape $r,$ where $L$ denotes the line through 
the origin and the endpoint of the SRW after $n$ steps.
Moreover, the expectation  
 $\bE_{\sip, \cV, \cL_r} = \bE_{\bQ^{\sip, \cV, r}_n}$
relative to the radial part is calculated as in 
$(\ref{axnumbershalfr})$ followed by $(\ref{averageliner})$ 
with $\cL = \cL_r.$
\end{definition}

Let $\bE_{\sip, \cV, *(r)}$ denote expectation of the 
$d$-dimensional weakly self-avoiding cone process relative 
to $\cV$ in shape $r.$ In particular, we write 
$\bE_{\sip, \cV, *} = \bE_{\sip, \cV, *(1/2)}.$ 
Thus, the definition of this process depends on the choice
of $\cV$ and on $\Phi.$
Note that there is no unique such process since only 
the distribution of the radial component of the process is
prescribed and not even the distribution on the lines in $\cV$
is specified. Consequently, there will be several ways to choose the
set $\cV.$  Importantly though, the shape carries much information.


\section{Expected Distances}
\setcounter{equation}{0}  

We turn to a technical lemma that engages a condition  
and a couple more definitions. 
The main players in this condition are bounded numbers $a_x$ 
that depend on $x>0$ and will be substituted by the
numbers $a_x(\cL_{r}),$ especially,  $a_x(\cL_{1/2}),$  
shortly. Recall that the latter are bounded in $x$ and
that we assumed that there is some number $\zeta >0,$
independent of $\sip,$ so that $\sip a_x(\cL_{r}) \geq \zeta$ 
for every $0 \leq x \leq n.$ 
Define
\begin{eqnarray}
             \label{mux}
   \mu_x & = & (\sip a_x)^{1/2} \, n^{3/4}  \\
             \label{functionk}
    q(x) & = &  \exp \{ - \sip \frac{a_x}{2}  n^{1/2} \} 
\end{eqnarray}   
for every $n \geq 0,$ $\sip > 0,$ and $x $ in $[0,n].$
Since $a_x $ is bounded in $x,$
for suitably small $\varepsilon \geq 0$ and for $ \gamma>0,$ 
we may define 
\begin{eqnarray} 
  \label{regioninter}
      r_1 & = & r_1(\varepsilon, \gamma) = 
           \sup \{ x \in [0, n] : 
          x \leq \gamma \mu_x n^{- \varepsilon} \} \nonumber \\
      r_2 & = & r_2( \gamma)  = 
           \sup \{ x \in [0, n] :   x \leq \gamma \mu_x  \} . 
\end{eqnarray} 
Thus, $r_2( \gamma) = r_1(0, \gamma).$

\smallskip
{\sc {\bf Condition D.} }
For any suitably small $\varepsilon \geq 0,$ 
there exist some $\gamma >0$ and $ \rho_* > 0$ such that 
\begin{equation}
  \label{conditiond}
   \int_{r_1}^{n} \, x \, q(x) \,  d \bP_{\chi_n}(x)
    = \rho_n \, \int_0^{r_1} \, x \, q(x) \,  d \bP_{\chi_n}(x) 
\end{equation}
with $\rho_n \geq \rho_* $ for all sufficiently large $n.$

\medskip
Note that if $ \int_{0}^{r_2} \, x \, q(x) \,  d \bP_{\chi_n}(x)
= o( \int_{r_2}^{n} \, x \, q(x) \,  d \bP_{\chi_n}(x))$
as $n \rightarrow \infty,$ then $\varepsilon =0$ and $\rho_n
\rightarrow \infty.$
In addition, observe that, in light of the expression in 
(\ref{functionk}) for $q(x),$ \mbox{Condition D} guarantees 
that $a_x$ not be constant in $x$ and $\sip >0.$ 
Throughout the paper, we shall be careful about whether 
constants in $n$ and/or $x$ depend on $\sip$ or not and  
indicate this. 

In the next result, drawn from {\sc Hueter} \cite{hue1},
the $a_x$ are some general numbers that obey the
stated assumptions.

\begin{lemma}[\cite{hue1}, Lemma 1]
{\bf (Exponent of Expected Radial Distance equals $3/4$)}
    \label{twointegrals}
Let $\sip > 0$ and $d \geq 1.$
Assume that the $a_x$ are bounded 
numbers that depend on $x,$ are such that there 
is some number $\zeta >0$ so that $\sip a_x \geq \zeta$ 
for every $0 \leq x \leq n,$ and that satisfy Condition D 
in $(\ref{conditiond})$ for some $\varepsilon \geq 0$ and $\gamma >0.$
Define
\begin{eqnarray}
    \label{threeintegrals}
 I_n  & = &  \int_0^n \, x \, q(x) \,  d \bP_{\chi_n}(x)  \\
 g(n) & = &  \int_0^n \, (a_x)^{1/2} q(x) \,  d \bP_{\chi_n}(x),
                              \nonumber 
\end{eqnarray}
where $q(x)$ is defined in $(\ref{functionk}).$
Then there are some constants $ M < \infty$ and $c(\rho_*) > 0$ 
(both independent of $\sip$) such that as $n \rightarrow \infty,$ 
\begin{equation}
  \label{inint}
  \gamma  \, c(\rho_*) \, \sip^{1/2} \, 
   n^{3/4 - \varepsilon} \, (1+o(1))
  \leq \frac{I_n}{g(n)} \leq  M \, \sip^{1/2} \, 
   n^{3/4} \, (1+o(1)). 
\end{equation}
\end{lemma}

\begin{proof}
We do not reproduce the proof that is presented in
{\sc Hueter} \cite{hue1}, Lemma 1, in two dimensions.
\end{proof}

The next result is as well borrowed from {\sc Hueter}
\cite{hue1}.

\begin{lemma}[\cite{hue1}, Lemma 2]
{\bf (The $a_x(\cL_{1/2})$ satisfy Condition D)} 
   \label{conditiondsatisfied}
Let $d > 1.$ 
If $\Phi$ has circular shape for sufficiently large $n,$
then the $a_x(\cL_{1/2}),$ 
defined in $(\ref{axnumbershalfr})$ when $\cL = \cL_{1/2}$ 
and $r=1/2$ satisfy Condition D in $(\ref{conditiond})$
for $\varepsilon = 0$ and $\gamma >0,$ independent of $\sip$
as $\sip \rightarrow \infty.$
\end{lemma}

\begin{proof}
The idea of proof is to violate Condition D 
and to take this assumption to a contradiction to the 
one that $\Phi$ be circular.
The proof is omitted here and can be found in {\sc Hueter}
\cite{hue1}, Lemma 2. It runs in parallel with the proof 
of Lemma \ref{conditiononed} that we will present in Section 6.
\end{proof}

\begin{proposition}[\cite{hue1}, Proposition 3]
{\bf (Expected Distance Along Cones with Order 
$n^{1/2}$ SILT)} 
    \label{distmainterm}  
Let $d >1$ and $\sip > 0.$ There are some constants
$ 0 < \gamma_* \leq  M < \infty$ (independent of $\sip$ as
$\sip \rightarrow \infty$ and $M$ independent of $\sip >0$ as well) 
such that as $n \rightarrow \infty,$   
\begin{eqnarray*}
  \bE_0 ( \chi_n \, e^{- \sip \philh})
          & = & K(n) \, n^{3/4} \sip^{1/2} g(n)(1+o(1))
\end{eqnarray*} 
for $ \gamma_* \leq K(n) \leq M,$
where $g(n)$ was defined in $(\ref{threeintegrals}).$
\end{proposition}

\begin{proof}
Since the proof is rather short, we present it here again.
Observe that, in view of Lemma \ref{conditiondsatisfied}, 
the $a_x = a_x(\cL_{1/2})$ satisfy Condition D in (\ref{conditiond})
for $\varepsilon = 0$ and $\gamma >0.$ 
Combining the observations preceding 
(\ref{distaverpenweight}) together with 
(\ref{axnumbershalfr}) and (\ref{averageliner}) 
and Lemma \ref{twointegrals}
leads to, as $n \rightarrow \infty,$   
\begin{eqnarray}
  \bE_0 ( \chi_n \, e^{- \sip \philh } )
           & = &  \bE_0 (\chi_n \,  \bE_{\Phi \vert \chi_n} ( \,
                       \vert \cL_{1/2}(\Phi) \vert^{-1}  
                          \sum_{L \in \cL_{1/2}(\Phi)} 
                     e^{- \sip \vert \cC_L \vert  } \vert 
                         \chi_n =x )) \nonumber \\*[0.2cm]
               & = &  \label{cesaroaverage}
                   \int_0^n \, x \, \bE_{\Phi \vert \chi_n} ( \,
                       \vert \cL_{1/2}(\Phi) \vert^{-1}  
                          \sum_{L \in \cL_{1/2}(\Phi)} 
                     e^{- \sip \vert \cC_L \vert  } \vert 
                         \chi_n =x ) \, d \bP_{\chi_n}(x)
                         \nonumber \\*[0.2cm]
               & = &  
                      \int_0^n \, x  \, (  
                     \, 
                      \int_{\bZ^d} \, 
                      \vert \cL_{1/2}(\varphi) \vert^{-1}  
                        \sum_{L \in \cL_{1/2}(\varphi)} 
                  e^{- \sip \vert \cC_L \vert  } \, 
                         d \bP_{\Phi \vert \chi_n}(\varphi \vert x) )
                           \, d \bP_{\chi_n}(x)
                        \nonumber    \\*[0.2cm]
               & = &  \label{evalpalmone}  
                   \int_0^n \, x  \, 
                            \exp \{ - \sip a_x n^{1/2}/2 \}
                         \, d \bP_{\chi_n}(x)
                          \\*[0.2cm]
                  & = &  \label{gnintegral} 
                   \int_0^n \, x  \, q(x) \, d \bP_{\chi_n}(x)
                          \\*[0.2cm]
                   & = &  K(n) n^{3/4} \sip^{1/2} g(n)(1+o(1))  
\end{eqnarray} 
for $ \gamma_* \leq K(n) \leq M,$ 
where to obtain the last two lines of the display, we apply 
\mbox{Lemma \ref{twointegrals},} 
with $\gamma_* = \gamma c(\rho_*), $ $\varepsilon=0,$ and with
the $a_x$ being bounded and such that there 
is some number $\zeta >0$ so that $\sip a_x \geq \zeta$ 
for every $0 \leq x \leq n.$ These two properties of $a_x$
may be seen as follows. First, since, by (\ref{halfline}),
$ 2 \vert \cC_L \vert / n^{1/2}$ is in $  [a_1, a_2 ],$
the average of the exponential 
terms $ \exp \{- \sip \vert \cC_L \vert \}$ over all lines
in $\cL_{1/2}(\varphi)$ may be rewritten as 
$\exp \{ - \sip a_x $ $ n^{1/2}/2 \},$ say,
for some number $a_x \in  [a_1, a_2 ],$ 
depending on $x.$ 
In particular, the $a_x$ are bounded. 
Additionally, we assumed (remark following display (\ref{halfline}))
that $a_1 \sip$ is a positive number independent of $\sip,$
thus, there is some number $\zeta >0$ so that $\sip a_x \geq \zeta$
for all $x.$ 
This completes our proof.
\end{proof}

We collect two propositions and key ingredients 
to our main results.

\begin{proposition}{\bf (Upper Bound for $ \bE_{\sip} \chi_n$)} 
   \label{msdupper}
Let $\sip >0$ and $d>1.$ There is some constant 
$ M_*(d) = M_*(d, \sip) < \infty$ (made precise below) such that 
as $n \rightarrow \infty,$  
\begin{eqnarray*} 
  \bE_{\sip} (\chi_n) & \leq &  M_*(d) \, (1 + o(1)) \,
              \max(n^{1/4+1/d}, n^{1/2}),
\end{eqnarray*}
where $ M_*(d)$ is uniform in $\sip$ for $d = 2$ and $d \geq 5$ and 
may depend on $\sip$ for $d =3,4.$ 
\end{proposition}
 
\begin{proof}
It will suffice to prove that, for $\cV$ in circular shape,
as $n \rightarrow \infty,$
\\*[0.1cm]
\mbox{} \qquad 
(I) \, \, $ \bE_{\sip, \cV, *} (\chi_n) 
           \leq M_*(d) \, (1+o(1)) \, \max (n^{1/4 + 1/d}, n^{1/2})$
  for $ M_*(d) < \infty$ and
\\*[0.05cm] \mbox{} \qquad 
(II) \,  $  \bE_{\sip}(\chi_n)  \leq  \bE_{\sip, \cV, *}(\chi_n)
                                       \, (1 + o(1)).$

\smallskip
{\bf Part (I).}
Assume that $\cV$ is $1/2$-shaped for all sufficiently large $n.$ 
First, we fix the size of $\cV.$
Since, up to scaling by a factor between $1$ and $\sqrt{d},$
any direction in ${\bf R}^d$ is the same for the
weakly self-avoiding cone process relative to $\Phi,$
we choose the lines in $\cV$ uniformly distributed over some
$d$-sphere such that there are of order $n^{1/d}$ lines along 
each side of the smallest $d$-cube that contains the $d$-sphere.
In other words, $ \vert \cV \vert$ is of order
$n^{(d-1)/d} = n^{1-1/d},$ say,
$ \vert \cV \vert = v_n \, n^{1-1/d}$ for $v_1(\sip) \leq v_n \leq v_2,$ 
for all sufficiently large $n,$ where, in view of (\ref{boundsonlr}),
we can choose the two constants $ 0 < v_1= v_1(\sip) 
\leq v_2 < \infty$ so that
$v_2$ is independent of $n$ and $\sip$ for each $\sip >0$
and $v_1$ is independent of $n$ but may depend on $\sip,$
even as $\sip \rightarrow \infty.$  

Next, a consequence of the arguments in the proof of Proposition 4
in {\sc Hueter} \cite{hue1} is that as $n \rightarrow \infty,$
\begin{eqnarray}
   \label{competingterms}
      \bE_{\sip, \cV, *}(\chi_n)  & \leq & 
      \bP_{\Phi} ( L \in \cL_{1/2} )
              \,  \frac{   \max_{\cL \subset \cL_{1/2} }  \,
             \bE_0 ( \chi_n \, e^{- \sip J_n^{\cL}} )}
           {\bE_0 ( e^{- \sip  J_n^{\cL_{1/2}}}  )}\, (1 + o(1))
                \nonumber \\*[0.05cm]
              & & \mbox{} \, \,  \, + \, 
           \bP_{\Phi} ( L \in  \cL_{\emptyset} )
              \,  \frac{ \bE_0 ( \chi_n \, e^{- \sip 
                J_n^{\cL_{\emptyset}}} )}
           {\bE_0 ( e^{- \sip  J_n^{ \cL_{\emptyset} }}  )} \, ,
\end{eqnarray}
where $ \cL_{\emptyset}$ and $\cL_{1/2}$ are defined in (\ref{halfline}) 
and $ \bE_0 ( \chi_n \, e^{- \sip J_n^{\cL}} )$ 
is to be understood in the sense of definitions (\ref{averageliner})
and (\ref{distaverpenweight}).
Furthermore, we have seen in {\sc Hueter} \cite{hue1}, 
proof of \mbox{Proposition 4,} that as $n \rightarrow \infty,$
\begin{equation}
 \label{halfquotient}
               \frac{   \max_{\cL \subset \cL_{1/2} }  \,
             \bE_0 ( \chi_n \, e^{- \sip J_n^{\cL}} )}
           {\bE_0 ( e^{- \sip  J_n^{\cL_{1/2}}}  )}
       \leq M \, (1 + o(1)) \, ( \sip a_2)^{1/2} \, n^{3/4} ,
\end{equation}
where $\sip a_2 $ and $M$ are finite constants that
do not depend on $\sip$ (uniform in $\sip$). 
Let us devote a moment to bound the last term 
in (\ref{competingterms}).
A routine exercise yields that as $n \rightarrow \infty,$
\begin{equation}
  \label{firstsrw}
   \frac{  \bE_0 ( \chi_n \, e^{- \sip J_n^{\cL_{\emptyset}}} )}
           {\bE_0 ( e^{- \sip  J_n^{\cL_{\emptyset}}}  )} 
       = (\frac{2}{\pi})^{1/2} \, n^{1/2} (1 + o(1)).
\end{equation}
The probability  $ \bP_{\Phi} ( L \in \cL)$ may be interpreted
as a Palm probability, that is,
\begin{equation}
   \label{palmlinesr}   
     \bP_{\Phi} ( L \in \cL)  =
           \frac{ \bE_{\Phi} \sum_{ L \in \cV} 1_{\cL}(L) }
                 { \vert \cV \vert }
          = \frac{ \bE_{\Phi} \vert \cL \vert}
                 { \vert \cV \vert } .
\end{equation} 
Therefore, since $ v_1 \, n^{1- 1/d} \leq \vert \cV \vert,$
by virtue of (\ref{boundsonlr}), 
\begin{equation}
  \label{halfprobup}
  \bP_{\Phi} ( L \in  \cL_{1/2}) \leq \frac{2 b_2}{ a_1 v_1} 
        \, n^{1/d -1/2} ,
\end{equation}
where $2 b_2 / (a_1 v_1) $ may depend on $\sip$ since $v_1$ does. 
Hence, 
(\ref{competingterms}),
(\ref{halfquotient}), 
(\ref{firstsrw}), and
(\ref{halfprobup}) can be summarized as
\begin{eqnarray}
   \label{distancerepresent}
 \bE_{\sip, \cV, *} (\chi_n)  
           & \leq &     
               \label{twoterms} 
             M_*(d) \, (1+o(1)) \, \max (n^{1/4 + 1/d}, n^{1/2}) 
\end{eqnarray} 
as $n \rightarrow \infty$
for $ M_*(d) = M (\sip a_2)^{1/2} (2 b_2 /(a_1 v_1)) +  (2/\pi)^{1/2}
< \infty$ for $ d =3,4,$ 
 $M_*(d) = M (\sip a_2)^{1/2}$ for $d=2,$
and $ M_*(d) = (2/\pi)^{1/2}$ for $d \geq 5.$ In summary,
$M_*(d)$ is uniform in $\sip$ for $d=2$ and $d \geq 5$
but may depend on $\sip,$ even as $\sip \rightarrow \infty,$
for $d =3,4.$ 
This completes the verification of (I) along with
the asymptotic evaluation of its righthand side.

\smallskip
{\bf Part (II).}
Demonstrating (II) will finish our proof. 
This portion is much as given in {\sc Hueter} \cite{hue1},
part (II) of the proof of Proposition 4. We sketch an
outline.
Recall that $\bE_{\sip, \cV, *(r)}$ denotes 
expectation of the $d$-dimensional 
weakly self-avoiding cone process relative 
to $\cV$ in shape $r.$ 
In order to compare $\bE_{\sip, \cV, *(r)}(\chi_n)$ and 
$\bE_{\sip}(\chi_n),$ the strategy will be to show that, for fixed
$J_n  \in [ b_1 n,  b_2 n ],$
the number of SRW-paths with $J_n $ whose point process $\Phi$ 
is $r$-shaped is larger than the number of SRW-paths with $J_n $ 
whose point process $\Phi$ is $s$-shaped (but not $r$-shaped) 
for $1/2 \leq r <s.$
We will continue to show that $\cL_s$ for $0 \leq s < 1/2$
plays a negligible role as well. 
In other words, most SRW-paths that
satisfy (\ref{rangejn}) arise from a $\Phi$ that is $1/2$-shaped.
Finally, we shall compare the centers of mass of the
weakly self-avoiding cone process and the weakly SAW.

{\bf (a) $\Phi$ prefers circular shape.}
Fix $J_n$ (and assume that $J_n \in [ b_1 n,  b_2 n ] $).
Partition the interval $[1/2,1]$ into $R$ subintervals of equal
length, that is, let $ 1/2 = r_0 < r_1 < r_2 < \ldots < r_R = 1.$
We are interested in comparing
the number of SRW-paths with $J_n$ whose point process $\Phi$ 
is $r_{k-1}$-shaped to the number of SRW-paths with $J_n$ whose 
point process $\Phi$ is $r_k$-shaped
(but not $r_{k-1}$-shaped). For this purpose,
we shall give an inductive argument over $k.$
Pick a SRW-path $\gamma$ of length $n$ with $J_n$ whose
point process $\Phi$ has shape $r_k.$ 
We will show that 
\mbox{(i) associated} with $\gamma,$ there is a large set $F_{\gamma}$ of 
SRW-paths whose realizations of $\Phi$ have shape $r_{k-1},$ and
(ii) two sets $F_{\gamma}$ and $F_{\gamma'}$ are disjoint for 
$\gamma \not = \gamma'.$ 
To see this, we {\em cut and paste} the path $\gamma$ as follows.
Let $P_{\gamma}$ denote the smallest parallelepiped that contains
the path $\gamma$ and let $l_{\gamma}$ denote the largest
integer less than or equal to the length of the longest side of
$P_{\gamma}.$  Divide $P_{\gamma}$
into sub-parallelepipeds whose sides are parallel to the sides of
$P_{\gamma}$ by partitioning the two longest sides of $P_{\gamma}$
into $n_f$ subintervals in the same fashion whose endpoints are
vertices of the integer lattice and by connecting the
two endpoints of the subintervals that are opposite to each
other on the two sides. Shift each of the sub-parallelepipeds
including the SRW-subpaths contained
by a definite amount between $1$ and $K$ ($K$: some constant) 
along one of the directions
of the shorter sides of $P_{\gamma}$ and reconnect the SRW-subpaths
where they were disconnected. In doing this, the shifts are
chosen such that the new path $\gamma'$ will have shape
$r_{k-1}$ and the total number of connections needed to reconnect
those subpaths equals a number $C_n$ that is constant in $k.$
Observe that such a choice of shifts exists.
When walking through the new path $\tilde{\gamma},$ 
because of the necessary extra steps to reconnect the subpaths,
the last several steps of $\gamma$ will be ignored. 
Note that this latter number of steps is independent of $k.$
Hence, if the pieces
to reconnect are self-avoiding, then $J_n$ is no larger after 
this cut-and-paste procedure than before.
This is always possible for otherwise we shift apart the 
sub-parallelepipeds such that they are sufficiently separated 
from each other.
Now, either we choose the reconnecting pieces such that
$J_n$ is preserved or we ``shift back'' (along the 
direction of the long sides of the parallelepiped) some or all of
the sub-parallelepipeds so that any two parallelepipeds 
overlap sufficiently to preserve $J_n$ and then 
reconnect the SRW-subpaths where they were disconnected. 
Again, we shift in such a fashion that the
total number of connections needed to reconnect
the subpaths equals $C_n.$ The number of these newly
constructed paths in $F_{\gamma}$ grows at least at the 
order that the number of ways does to choose $n_f$ locations 
(to shift) among $l_{\gamma}$ sites, which is a number larger
than $1$ for all large enough $n.$
Hence, the number of SRW-paths with $J_n$ whose 
point process $\Phi$ is $r_{k-1}$-shaped is larger
than the number of SRW-paths with $J_n$ whose 
point process $\Phi$ is $r_k$-shaped. Since
this argument can be made for every $1 \leq k \leq R$
and the number of SRW-paths with $J_n$ whose 
point process $\Phi$ has shape $r_R=1$ is at least $1,$
it follows that the number of SRW-paths with $J_n$ whose
point process $\Phi$ is $r$-shaped is maximal for $r=1/2.$ 

{\bf (b) It suffices to consider shapes $r$ with $r \geq 1/2.$}
This passage is identical to part (II)(b) in the proof of 
Proposition 4 in {\sc Hueter} \cite{hue1} and omitted here.

The considerations in (a) above also imply that both probability
distributions decay exponentially fast around their centers of mass. 
Combining this observation with the fact that
the shape of $\Phi$ relates the SILT of the weakly SAW to the one of
the weakly self-avoiding cone process provides that the two probability
distributions asymptotically have the same centers of mass (up to error
terms).
Together with these, the upshot of above passages (a) and (b) is that,
in comparing $ \bE_{\sip}(\chi_n)$ to $\bE_{\sip, \cV, *(r)}(\chi_n)$
for $ 0 \leq r \leq 1,$ it is enough to choose $r=1/2$ and to 
study the expected distance of the weakly self-avoiding cone 
process relative to $\Phi$ when in circular shape. 
Hence, in particular, we are led to
$$
   \bE_{\sip}(\chi_n)  \leq  \bE_{\sip, \cV, *}(\chi_n) (1 + o(1))
$$
as $n \rightarrow \infty.$
This accomplishes the proof of (II), and thus, ends the proof.
\end{proof}
 
\begin{proposition}{\bf (Lower Bound for $ \bE_{\sip} \chi_n$)} 
   \label{msdlower}
Let $\sip >0$ and $d>1.$ There is a constant $m(d)=m(d, \sip) >0$ 
(made precise below) such that as $n \rightarrow \infty,$  
\begin{eqnarray*}  
         \bE_{\sip} (\chi_n)  & \geq & 
               m(d)  \, (1 + o(1)) \,
              \max(n^{1/4+1/d}, n^{1/2}) ,
\end{eqnarray*}
where $ m(d)$ is uniform in $\sip$ for $d \geq 5$ and 
may depend on $\sip$ for $d =2, 3, 4.$ 
\end{proposition}  

\begin{proof}   
Fix $\rho >0$ and $\delta >0.$
From the reasoning in {\sc Hueter} \cite{hue1}, proof of
Proposition 5, and part (II)(b) in the proof of Proposition 
\ref{msdupper}, we collect 
\begin{equation}
  \label{lowercomp}
      \bE_{\sip, \cV, *}(\chi_n) \leq \bE_{\sip}(\chi_n) 
\end{equation}
and that there are two competing terms involved to bound 
$\bE_{\sip, \cV, *}(\chi_n)$ from below.
The dimension $d$ decides which one is maximal.
(Hence, for the rest of the proof, we may assume that 
$\Phi$ has circular shape.)
Namely, we have, as $n \rightarrow \infty,$
\begin{eqnarray}
    \label{minpalm}
\bE_{\sip, \cV, *}(\chi_n) 
               & \geq &   (1 + o(1)) \,
                  \max [\, \bP_{\Phi} ( L \in \cL_{1/2} )
               \,  \frac{   \min_{\cL \subset \cL_{1/2}}  \,
             \bE_0 ( \chi_n \, e^{- \sip J_n^{\cL}} )}
           {\bE_0 ( e^{- \sip  J_n^{\cL_{1/2}}}  )}  \, ,
         \\*[0.1cm] & & \mbox{} \, \, \, \mbox{} 
       \bP_{\Phi} ( L \in \cL_- \cup \cL_{\emptyset}) \, 
      \frac{  \min_{\cL \subset \cL_- \cup \cL_{\emptyset}}  \, 
         \bE_0 ( \chi_n \, e^{- \sip J_n^{\cL}} )}
        {\bE_0 ( e^{- \sip  J_n^{\cL_- \cup \cL_{\emptyset}}}  )} \, ]   
      \nonumber \\*[0.2cm]  & \geq & (1 + o(1)) \,
                  \max [\, \bP_{\Phi} ( L \in \cL_{1/2} )
               \,  \frac{   \min_{\cL \subset \cL_{1/2}}  \,
             \bE_0 ( \chi_n \, e^{- \sip J_n^{\cL}} )}
           {\bE_0 ( e^{- \sip  J_n^{\cL_{1/2}}}  )}  \, ,
        \nonumber \\*[0.1cm] & & \mbox{} \, \, \, \mbox{} 
         \bP_{\Phi} ( L \in \cL_- \cup \cL_{\emptyset}) \,
        \frac{ \bE_0 ( \chi_n \, e^{- \sip J_n^{\cL_{\emptyset}}} )}
           {\bE_0 ( e^{- \sip  J_n^{ \cL_{\emptyset}}}  )} \, ] ,
        \nonumber 
\end{eqnarray}
where $ \cL_{\emptyset}$ and $\cL_{-}$ are defined in (\ref{halfline})
and the minimum $\min_{\cL \subset \cL_{1/2}}$
is over subsets $\cL \subset \cL_{1/2}$ 
that form a subset of $\cV$ that is circular for sufficiently 
large $n.$ 
Moreover, in {\sc Hueter} \cite{hue1}, proof of
Proposition 4, we arrived at 
\begin{equation}
  \label{halfquotlower}
  \frac{   \min_{\cL \subset \cL_{1/2}}  \,
             \bE_0 ( \chi_n \, e^{- \sip J_n^{\cL}} )}
           {\bE_0 ( e^{- \sip  J_n^{\cL_{1/2}}}  )}
         \geq (1 + o(1)) \, \gamma_* \, ( \sip a_1)^{1/2}
            \, n^{3/4} 
\end{equation}
as $n \rightarrow \infty,$ where $\gamma_* $ is independent
of $\sip$ as $\sip \rightarrow \infty$ and 
$\sip a_1 $ is a positive number that is independent
of all $\sip >0.$
The fact that we assumed $\Phi$ to be $1/2$-shaped
together with (\ref{boundsonlr}) for $r=1/2$ yields 
$ \vert \cL_{1/2} \vert \geq (b_1/a_2) \, n^{1/2 - \rho}.$
As we pointed out at the outset of the proof of Proposition
\ref{msdupper}, part (I), we can choose
$ \vert \cV \vert = v_n \, n^{1-1/d},$ where $v_n \leq v_2$ 
for all sufficiently large $n$ and
$v_2$ is independent of $\sip.$ 
Consequently, we end up with 
\begin{equation}
   \label{halfprob}  
     \bP_{\Phi} ( L \in \cL_{1/2})  =
           \frac{ \bE_{\Phi} \sum_{ L \in \cV} 1_{\cL_{1/2}}(L) }
                 { \vert \cV \vert } > m_* \, n^{1/d - 1/2 - \rho}
\end{equation}
for $m_* = b_1/(v_2 a_2) >0 $ and every sufficiently large $n.$ 
Note that $m_*$ may depend on $\sip,$ even as $\sip
\rightarrow \infty,$ since $b_1/a_2$ does.

Furthermore, recall $ \cL_{+}$ and $\cL_{1/2 \pm}$ 
from (\ref{halfline}). 
It is obvious that $ \vert \cL_{-} \cup \cL_{\emptyset} \vert
= \vert \cV \vert - $ $ \vert  \cL_{1/2 \pm} \cup \cL_{+} \vert.$ We
claim that $\vert  \cL_{1/2 \pm} \cup \cL_{+} \vert 
\leq (R+1) (2 b_2 / a_1) \, n^{1/2} = R' \, n^{1/2},$ where 
$R' = (R+1) (2 b_2 / a_1)$
is a finite constant, independent of $n.$ 
To see this, fix some integer $R >0$ and
partition the interval $[1/2 - 1/2R, 1],$ into $R+1$ 
subintervals of equal lengths, that is, 
$1/2-1/2R =r_0 < r_1 < \ldots < r_R =1-1/2R < r_{R+1} =1$
with $r_{k+1} = r_k + 1/2R$ for all $k.$ Consider the sets 
$\cL_{r \star} =
 \{L \in \cV: \, 2 \vert \cC_L \vert
      \in [ a_1 n^r,  a_2 n^{r + \delta}] \} $ for $\delta= 1/2R$
and $r=r_0, r_1, \ldots, r_{R}.$ 
We have $  \cL_{1/2 \pm} \cup \cL_{+} \subset \cup_{k=0}^R \cL_{r_k \star}.$
Therefore, in view of (\ref{boundsonlr}) for $r \geq 1/2- 1/2R,$  
we find $\vert \cL_{1/2 \pm} \cup \cL_{+} \vert \leq
(R+1) (2b_2/a_1) n^{1/2 + 1/2R}.$

Hence, when $d \geq 3,$ we can choose $R$ suitably large such 
that we obtain 
$\vert \cL_{1/2 \pm} \cup \cL_{+} \vert  = o( \vert \cV \vert )$ 
as $n \rightarrow \infty$ (since $\vert \cV \vert = v_n \, n^{1-1/d}$). 
Whence, when $d \geq 3,$ as $n \rightarrow \infty,$
\begin{eqnarray}
  \label{lowsiltlines}
  \bP_{\Phi} ( L \in \cL_- \cup \cL_{\emptyset}) 
      & = & \frac{ \vert \cV \vert - 
         \vert  \cL_{1/2 \pm} \cup \cL_{+} \vert 
         }{ \vert \cV \vert }  \geq (1 - o(1)) .
\end{eqnarray} 
Finally, in light of 
(\ref{firstsrw}),
(\ref{minpalm}),
(\ref{halfquotlower}), 
(\ref{halfprob}), and
(\ref{lowsiltlines}), 
we obtain, as $n \rightarrow \infty,$
\begin{eqnarray}
    \label{lbpalmone} 
 \bE_{\sip} (\chi_n) & \geq & 
               m(d) \,  (1 + o(1)) \, \max (n^{1/4+1/d - \rho}, n^{1/2})
\end{eqnarray}  
for $m(d) =\max( m_* \gamma_* (\sip a_1)^{1/2} , (2/\pi)^{1/2} )$
for $d=2,3,4$ and 
$m(d) =(2/\pi)^{1/2}$ for $d \geq 5.$
We summarize to say that $m(d)$ is uniform in $\sip$ for $d \geq 5$
but may depend on $\sip,$ even as $\sip \rightarrow \infty,$
for $d =2,3,4.$ 
Since $\rho > 0$ was arbitrary, the announced lower bound 
for $ \bE_{\sip} (\chi_n) $ is an immediate consequence.
\end{proof}


\section{Distance Exponents of the Self-Avoiding Walk}
\setcounter{equation}{0} 

We recall the numbers $ \mu (d) = \max (1/4+1/d , 1/2)$ 
for every integer $d >1$ from (\ref{valueofexponent}).
A number of arguments towards uniform bounds in 
$\sip$ as $\sip \rightarrow \infty$
together with Propositions \ref{msdupper} and \ref{msdlower}
will demonstrate that the values of the
distance exponents extend to the SAW.

\begin{corollary}
  \label{liminfsupmsd}
Let $\sip > 0$ and $d>1.$ 
There are some constants $0 < \rho_1(d) = \rho_1(d, \sip)
\leq \rho_2(d) = \rho_2(d, \sip) < \infty$  such that 
$$  
   \rho_1(d)  \leq
     \liminf_{n \rightarrow \infty}
         n^{- \mu(d)}  \, \bE_{\sip} (\chi_n) 
        \leq  \limsup_{n \rightarrow \infty}
        n^{- \mu(d)}  \, \bE_{\sip} (\chi_n) 
         \leq \rho_2(d),
$$ 
where $\rho_1(d)$ is uniform in $\sip$ for $d \geq 5$ and 
may depend on $\sip$ for $d \leq 4$ and
$\rho_2(d)$ is uniform in $\sip$ for $d = 2$ and $d \geq 5$ and 
may depend on $\sip$ for $d =3,4.$
In particular, the self-avoiding walk in $\Zd$ for $d \geq 2$
has distance exponent $ \max (1/4+1/d , 1/2).$
\end{corollary}

\begin{proof}
The statements on the weakly SAW are immediate consequences
of Propositions \ref{msdupper} and \ref{msdlower}.
The reasoning to verify that the SAW has the same distance
exponents in $\Zd$ for $d \geq 2$ as the weakly SAW is
identical to the one presented in {\sc Hueter} \cite{hue1},
proof of \mbox{Corollary 1.} The part to establish the lower bound $3/4$
for the distance exponent of the SAW was slightly more involved.
For $d=3,4,$ in fact, we apply those lines of arguments to both
the upper and lower bounds for the distance exponent of the SAW.
The details are omitted here.
\end{proof}

This completes the proof of Theorem \ref{sawdistance}. 
The next result accomplishes Theorem \ref{rmsd}.

\begin{corollary}
  \label{variance}
Let $\sip > 0$ and $d>1.$ 
There are some constants $0 < \rho_3(d) = \rho_3(d, \sip)
\leq \rho_4(d) = \rho_4(d, \sip) < \infty$  such that 
$$
   \rho_3(d) \leq
     \liminf_{n \rightarrow \infty}
       n^{- 2 \mu(d)} \, \bE_{\sip} (\chi_n^2) 
        \leq  \limsup_{n \rightarrow \infty}
       n^{- 2 \mu(d)} \, \bE_{\sip} (\chi_n^2) 
         \leq \rho_4(d),
$$ 
where $\rho_3(d)$ is uniform in $\sip$ for $d \geq 5$ and 
may depend on $\sip$ for $d \leq 4$ and
$\rho_4(d)$ is uniform in $\sip$ for $d = 2$ and $d \geq 5$ and 
may depend on $\sip$ for $d =3,4.$
In particular, the MSD exponent of the SAW in $\Zd$ for $d \geq 2$
equals $ \, 2 \max (1/4+1/d , 1/2).$
\end{corollary}
 
\begin{proof}
We shall argue in the setting of the weakly SAW and 
point out that parallel lines to the ones given before allow to
extend the results about the MSD exponent to the SAW.
First, the lower bound follows from the inequality
$\bE_{\sip} \chi_n^2 \geq (\bE_{\sip} \chi_n)^2$
and Corollary \ref{liminfsupmsd}. 

Hence, it is enough to verify the upper bound.
For this purpose, in light of (\ref{competingterms}),
it suffices to consider shapes $s$ of $\Phi$ for $s \geq 1/2$ 
in the case when $2 \leq d \leq 4,$ and, to collect
the contribution that would come from the second moment
of the distance of the SRW for $d \geq 5$ (compare also 
to (\ref{firstsrw})).
For $2 \leq d \leq 4,$ we remark the following.
It follows from refined considerations of those 
in Proposition \ref{msdupper},
part (II), that the number of SRW-paths that have
circular shape is exponentially larger (in $n$) than
the number of SRW-paths that have shape $s > 1/2.$
Also, the expected penalizing weight is exponentially
smaller for $\cL_s$ and $s > 1/2$
and decays exponentially fast in $n.$ Combining both of 
these observations implies that $\bQ_n^{\sip}$ decays 
exponentially fast in $n$ around $ \bE_{\sip} ( \chi_n ).$ 
Exponential decay holds for $d \geq 5,$ too.
In other words, for every $\epsilon >0,$ if we write 
$M_{\epsilon} = (\rho_2(d) \, n^{\mu(d)})^{1+ \epsilon}, $ we have
$$
  \bQ_n^{\sip} ( \chi_n > M_{\epsilon}) 
    \leq 
   \bQ_n^{\sip} ( \chi_n \leq M_{\epsilon}) 
         \, \exp \{  - \kappa(n, \epsilon) \}
$$
where $ \kappa(n, \epsilon) \rightarrow \infty $ as $n \rightarrow
\infty.$
Hence, since $\chi_n^2$ is bounded for every $n,$ we have
as $n \rightarrow \infty,$
\begin{eqnarray*}
   \bE_{\sip} ( \chi_n^2) & =  & 
                        (1 + o(1)) \, \sum_{k=0}^{M_{\epsilon}^2}  
                          \bQ_n^{\sip}( \chi_n^2 \geq k )  \\[0.1cm]
                & \leq & (1 + o(1)) \,  (M_{\epsilon}^2 +1) .
\end{eqnarray*}
Since $ \epsilon >0$ was arbitrary, we let $\rho_4(d) = \rho_2(d)^2$ to
wind up with the required result for \mbox{$d \geq 2.$}
\end{proof}

\begin{theorem} 
  \label{convexhull}
Let $\sip > 0$ and let $R_n$ denote the radius of the convex
hull of the SRW-path $S_0, S_1, \ldots, S_n$ 
in $\Zd$ for $d \geq 2.$ Then $R_n$ satisfies
all statements in Corollaries $\ref{liminfsupmsd}$
and $\ref{variance}$ with $\chi_n$ replaced by $R_n.$
\end{theorem}         
 
\begin{proof}
We are interested in the {\em maximal} distance
of $S_0, S_1, \ldots, S_n$ along any line
rather than the distance of the position of
$S_n$ from the starting point. 
The reflection principle gives the upper bound 
$d \bP_{R_n}(x)/dx \leq 2 \, d \bP_{\chi_n}(x)/dx$
whereas the lower bound $d \bP_{R_n}(x)/dx \geq  d \bP_{\chi_n}(x)/dx$ 
is readily apparent. The results now follow.
\end{proof}              

\medskip
{\bf Remark (Transition $\sip \rightarrow 0$).} 
For $d \leq 3,$ the transition $\sip \rightarrow 0$ looks dramatically
different from the transition $\sip \rightarrow \infty.$ 
The principal reasons are reflected upon the expressions in
(\ref{competingterms}) and (\ref{minpalm}).
Let us briefly glance at what happens 
when $\sip =0.$ In that fictive case (since the results were
proved under the assumption $\sip >0$),
both terms in (\ref{competingterms}) are of asymptotic order $n^{1/2},$
and so is the term in (\ref{minpalm}). Since this is drastically
different from the case $\sip >0,$ at least when $d=1,2,$ and $3,$ 
in which case the asymptotic order in $n$ of the largest term is 
$n^{\mu(d)} > \! \! > n^{1/2},$ we observe a discontinuity
of the expected distance measures and distance exponents 
of the weakly SAW as $\sip \rightarrow 0$ for $ d \leq 3.$
In contrast, the case $\sip \rightarrow \infty$ behaves as any
case for fixed $\sip.$

 
\section{One-Dimensional Weakly SAW}
\setcounter{equation}{0} 

This paragraph handles the case $d=1.$ We note in passing
that the one-dimensional SAW is not interesting.
Recall the numbers $a_x(\cL_1)$ from (\ref{axnumbershalfr})
with $r=1.$ They take values in $[a_1, a_2].$ For the remainder
of the paper, since $d=1,$ we assume that $a_1 = b_1 $ and 
$a_2 = b_2.$ Recall that $a_1 \sip$ depends on $\sip$ while
$a_2 \sip$ does not. Define
\begin{eqnarray}
   \label{muxone}
   \mu_x(1) & = &  (\sip a_x)^{1/2} \, n \\
   \label{functionkone}
        q_1(x) & = &  \exp \{ - \sip \frac{a_x}{2}  n \} 
\end{eqnarray}   
for every $n \geq 0,$ $\sip > 0,$ and $x $ in $[0,n].$
Similarly as in (\ref{regioninter}), 
for suitably small $\varepsilon \geq 0$ and for $ \gamma>0,$  
we may define
\begin{eqnarray} 
  \label{ronertwo}
      \hat{r}_1 & = & \hat{r}_1(\varepsilon, \gamma) = 
           \sup \{ x \in [0, n] : 
          x \leq \gamma \mu_x(1) n^{- \varepsilon} \} \nonumber \\
      \hat{r}_2 & = & \hat{r}_2(\gamma) = \hat{r}_1(0 , \gamma).
\end{eqnarray} 
Then the numbers $a_x(\cL_1),$ defined in (\ref{axnumbershalfr}) 
when $r=1,$ satisfy a hypothesis analogous to Condition D 
in (\ref{conditiond}).

\begin{lemma}
{\bf (Condition satisfied by $a_x(\cL_1)$)} 
   \label{conditiononed}
Let $d = 1.$ Then the $a_x(\cL_1)$ obey the following condition
for every $\varepsilon \geq 0:$

{\bf Condition $\tilde{\mbox{D}}$.}
For any suitably small $\varepsilon \geq 0,$ 
there exist some $\gamma >0$ and $ \omega_* > 0$ such that 
\begin{equation}
  \label{conditiondone}
   \int_{\hat{r}_1}^{n} \, x \, q_1(x) \,  d \bP_{\chi_n}(x)
    = \omega_n \, \int_0^{\hat{r}_1} \, x \, q_1(x) \,  d \bP_{\chi_n}(x) 
\end{equation}
with $\omega_n \geq \omega_* $ for all sufficiently large $n.$
\end{lemma}
 
\begin{proof}
Showing this statement requires no more than a small number of 
modifications to the proof of Lemma 2 in {\sc Hueter} \cite{hue1}
and can do without the assumption that $\Phi$ be circular. 
Fix some suitably small $\varepsilon >0.$
Let us invoke the notation that we introduced in the proof of
Lemma \ref{twointegrals}, that is, write
$ \bE_0 ( \chi_n \, e^{- \sip \philone } ) = I_n
=J_1(n) + J_2(n) + J_3(n),$ and in the same spirit,
$ \bE_0 ( e^{- \sip \philone } )  
= \tilde{J}_1(n) + \tilde{J}_2(n) + \tilde{J}_3(n),$
where the $r_i$ for $i=1,2,3$ are replaced by the $\hat{r}_i.$  
We need to show that there is some $\omega_* >0$ so that
$J_2(n) + J_3(n) = \omega_n J_1(n)$ with $\omega_n \geq \omega_* $ 
for all sufficiently large $n>0.$
We first show that $J_2(n) + J_3(n) \not = o(J_1(n))$
as $n \rightarrow \infty.$ 

For a moment, let us suppose in contrast
that $J_2(n) + J_3(n)  = o(J_1(n))$ as $n \rightarrow \infty$
so as to take this claim to a contradiction.
Thus,  $J_2(n)  = o(J_1(n))$ and $J_3(n)  = o(J_1(n))$ 
as $n \rightarrow \infty.$ It would
follow that $I_n = J_1(n)(1 + o(1))$ as $n \rightarrow \infty$
and
$\sum_{i=1}^3 \tilde{J_i}(n) = \tilde{J}_1(n) (1 + o(1)).$   
The probability measure $ \bQ^{\sip, \cV, 1}_n $
induces a one-dimensional process $W_n$ which has expectation
$\bE_{\sip, \cV, \cL_{1}}( \chi_n),$ call it 
$ \bE_{\sip, \cV, \cL_{1}}^{W}(\chi_n). $ 
Associate $W_n$ with the numbers $a_x(\cL_{1}).$ 

In view of the exponential form of the integrand of $I_n,$ 
our assumption would imply that there is a number
$z_n=z$ in $  [0, \hat{r}_1]$ that enjoys the property 
\begin{equation}
   \label{saddleone}
   \frac{ \bE_0 ( \chi_n \, e^{- \sip \philone}) }
        {\bE_0 ( e^{- \sip \philone } ) }  = 
    \frac{I_n}{ \bE_0 ( e^{- \sip \philone } ) }  = 
    (1 + o(1)) \, z 
\end{equation}
as $n \rightarrow \infty.$ In that event,
the function $a_x$ is minimal at $z=z_n,$ 
that is $a_z = \inf_{0 \leq y \leq \hat{r}_1} a_y$
for all sufficiently large $n.$
This may be seen as follows.
Define $k_1(x)   =  \exp \{ - (x^2 + \mu_x(1)^2 )/(2n) \},$   
let $a_0 >0$ and let $0 <\tau \leq a_0$ 
be some arbitrarily small number. If $a_{x_1} = a_0$ and
$a_{x_2} = a_0 - \tau \geq 0$ for $0 \leq x_1, x_2 \leq \hat{r}_1,$ 
then it follows that $k_1(x_1) < k_1(x_2)$ for all
sufficiently large $n.$

Now, for some suitably small $\tau = \tau(\sip) >0,$ define the set 
$$ 
  S_{\tau} = \{ x \in [0, \hat{r}_1] : a_x > a_z + \tau \} .
$$ 
Consider a modified process $\tilde{W}_n$ that is associated
with numbers $\tilde{a}_x$ with
$\tilde{a}_x = a_x$ for $x \in [0, \hat{r}_1] \setminus S_{\tau},$
$\tilde{a}_x = a_z + \tau$ for $x \in S_{\tau},$ and
$\tilde{a}_x = a_x + a(n)$ for $\hat{r}_1 <  x \leq n,$ where
$a(n)>0 $ is some suitable number, chosen so as to preserve the
distribution of $J_n.$ 
Thus, $\tilde{a}_x \leq a_z + \tau$ for $x \in [0, \hat{r}_1].$  
Observe that the modified process $\tilde{W}_n$ has the same
expectation $ \bE_{\sip, \cV, \cL_{1}}^{\tilde{W}}(\chi_n) =
\bE_{\sip, \cV, \cL_{1}}^{W}(\chi_n)$
as the process $W_n$ since, firstly, $q_1(x)$ in (\ref{functionkone}) 
was decreased on $(\hat{r}_1, n],$
and thus, $J_2(n) + J_3(n)$ and the corresponding part
$\tilde{J}_2(n) + \tilde{J}_3(n)$ of the integral
in the denominator of $\bE_{\sip, \cV, \cL_{1}}^{W}(\chi_n)$
were both decreased, and secondly, $J_1(n)$ is as before in 
view of (\ref{saddleone}). Note that adding 
a constant number of self-intersections to {\em all} 
realizations of this underlying weakly self-avoiding process 
$\tilde{W}_n$ does not change its probability measure.
Subtract the number $a_z$ from $\tilde{a}_x$ 
for every $0 \leq x \leq n,$ that is,
let $\hat{a}_x = \tilde{a}_x - a_z \geq 0 $ 
for every $0 \leq x \leq n.$
Thus, $\hat{a}_x \leq \tau$ for $x \in [0, \hat{r}_1]$
and  $\hat{a}_x$ is suitable on $(\hat{r}_1, n].$
The gotten process $\hat{W}_n$
associated with the numbers $\hat{a}_x$ has 
expectation  $ \bE_{\sip, \cV, \cL_{1}}^{\hat{W}}(\chi_n) =
\bE_{\sip, \cV, \cL_{1}}^{W}(\chi_n), $ too, the same
as do $W_n$ and $\tilde{W}_n.$
Since we may choose $\tau < b_1,$ we arrive at
$\hat{a}_x \leq \tau < b_1.$ But this contradicts
condition (\ref{rangejn}). We conclude that 
$J_2(n) + J_3(n) \not = o(J_1(n))$ as $n \rightarrow \infty.$ 
Since $\varepsilon >0$ was arbitrary,
it follows that, for every $\varepsilon >0,$
$J_2(n) + J_3(n) \not = o(J_1(n))$ as 
$n \rightarrow \infty.$

It remains to be shown that there is no subsequence $n_k$
such that $J_2(n_k) + J_3(n_k) = o(J_1(n_k))$
as $k \rightarrow \infty.$ From this it will follow that
there is some number $\omega_* >0$ that bounds $\omega_n$
from below with $n.$ But the same point can be made as
explained above when $n$ is replaced by $n_k$ everywhere.
Whence, we conclude that 
 (\ref{conditiondone}) is satisfied for every $\varepsilon >0.$
Hence, it must hold for $\varepsilon =0.$
Note that this implies that $\hat{r}_1 = \hat{r}_2$ 
and $J_2(n)= 0.$

Finally, we remark that $\gamma >0$ may be chosen 
uniformly over $\sip$ as $\sip \rightarrow \infty.$
This can be seen as follows.
Any of the asymptotic statements in a variant of Lemma  
\ref{twointegrals} and in the above lines of proof 
depend on expressions, for example, 
of the form $\sip^{1/2} n.$
Hence, if $N(\sip)$ is a threshold so that, for all
$n \geq N(\sip),$ a given expression in $n$ differs
from its corresponding limiting expression
by at most $\varepsilon$ (some fixed $\varepsilon$), it follows
that $N(\sip') \leq N(\sip)$ for $\sip < \sip'.$ As a
consequence of the fact that $\gamma >0$ may be chosen
uniformly in $n,$ the choice of $\gamma$ is uniformly
over $\sip >0$ as $\sip \rightarrow \infty$
(yet not as $\sip \rightarrow 0.$)
This accomplishes the proof. 
\end{proof}

\begin{theorem}
  \label{onedimension}
Let $\sip > 0$ and $d=1.$ Then there is a constant
$m_1(\sip) >0  $ that may depend on $\sip$
such that as $n \rightarrow \infty,$   
\begin{eqnarray*}
 \bE_{\sip} (\chi_n)  & = & (1 + o(1)) \, M_1(n) \, n 
\end{eqnarray*} 
for $ m_1(\sip) \leq  M_1(n) \leq 1.$ 
\end{theorem}

\begin{proof}
Clearly, it is sufficient to verify the lower bound for
$  \bE_{\sip} (\chi_n).$ 
Analogous lines to the one to obtain (\ref{minpalm}) lead to
\begin{eqnarray}
    \label{lowpalmone}
 \bE_{\sip}(\chi_n) & \geq &   
             (1 + o(1)) \, \bP_{\Phi} ( L \in \cL_1 )
       \,  \frac{ \bE_0 ( \chi_n \, e^{- \sip J_n^{\cL_1}} )}
           {\bE_0 ( e^{- \sip  J_n^{\cL_1}}  )} 
\end{eqnarray}  
as $n \rightarrow \infty.$ 
First, observe that $ \bP_{\Phi} ( L \in  \cL_{1}) \geq 1/2$ 
(there are only two half-lines that emanate from the origin in
${\bf Z}$).
Next, writing 
\begin{eqnarray*}
  \label{twointegralsone}
    I_n(1) & = &  \int_0^n \, x \, q_1(x) \,  d \bP_{\chi_n}(x) \\
    g_1(n) & = & \int_0^n \, (a_x)^{1/2}  q_1(x) \,  d \bP_{\chi_n}(x) 
\end{eqnarray*}  
and proceeding as to prove Lemma 1 in {\sc Hueter} \cite{hue1}
in connection with Lemma \ref{conditiononed} with $\varepsilon =0,$
we collect, as $n \rightarrow \infty,$   
\begin{equation}                                        
   \label{oneintegral}
     I_n(1) =  K_1(n) \, \sip^{1/2} \, g_1(n) \, n \, (1+ o(1))
\end{equation}
for $  0 < \gamma \, c(\omega_*) \leq K_1(n),$
where the positive constant $\gamma $ arises in 
Condition $\tilde{\mbox{D}}$ stated in Lemma \ref{conditiononed} 
and $c(\omega_*)$ is a positive constant, independent of $\sip.$
Both, $\gamma $ and $ c(\omega_*)$ may be chosen independently of $\sip $
as $\sip \rightarrow \infty,$ as we reasoned earlier.
In parallel to the arguments in proving 
\mbox{Proposition \ref{distmainterm}} together with the estimate
in (\ref{oneintegral}), we obtain, as $n \rightarrow \infty,$ 
\begin{eqnarray*}
  \bE_0 ( \chi_n \, e^{- \sip \philone } )
           & = &  \bE_0 (\chi_n \,  \bE_{\Phi \vert \chi_n} ( \,
                       \vert \cL_{1}(\Phi) \vert^{-1}  
                          \sum_{L \in \cL_{1}(\Phi)} 
                     e^{- \sip \vert \cC_L \vert  } \vert 
                         \chi_n =x )) \nonumber \\[0.2cm]
               & = &  
                   \int_0^n \, x \, \bE_{\Phi \vert \chi_n} ( \,
                       \vert \cL_{1}(\Phi) \vert^{-1}  
                          \sum_{L \in \cL_{1}(\Phi)} 
                     e^{- \sip \vert \cC_L \vert  } \vert 
                         \chi_n =x ) \, d \bP_{\chi_n}(x)
                         \nonumber \\[0.2cm]
               & = &  
                      \int_0^n \, x  \, ( \, 
                      \int_{\bZ} \, 
                      \vert \cL_{1}(\varphi) \vert^{-1}  
                        \sum_{L \in \cL_{1}(\varphi)} 
                  e^{- \sip \vert \cC_L \vert  } \, 
                         d \bP_{\Phi \vert \chi_n}(\varphi \vert x) )
                           \, d \bP_{\chi_n}(x)
                        \nonumber    \\[0.2cm]
               & = &  
                   \int_0^n \, x  \, 
                            \exp \{ - \sip a_x n /2 \}
                         \, d \bP_{\chi_n}(x)
                          \\*[0.2cm]
                  & = &   
                   \int_0^n \, x  \, q_1(x) \, d \bP_{\chi_n}(x)
                          \\*[0.2cm]
                   & = &  K_1(n) \, \sip^{1/2} \, g_1(n) \, n \, (1+o(1))  
\end{eqnarray*} 
for some $ 0 < \gamma_1 \leq K_1(n),$ where $\gamma_1 >0$ is 
independent of $\sip $ as $\sip \rightarrow \infty.$
Additionally, upon a similar but easier exercise, we gain
\begin{eqnarray*}
  \bE_0 ( e^{- \sip \philone } )
                            & = &  
                   \int_0^n   q_{1}(x) \, d \bP_{\chi_n}(x)
                          \\*[0.2cm]
                    & = &  h_{1}(n).
\end{eqnarray*}  
Note that $ g_1(n) /h_1(n) \geq (a_1)^{1/2} = (b_1)^{1/2}.$  
Hence, putting these two or three pieces together, along with
(\ref{lowpalmone}), yields as
$n \rightarrow \infty,$
\begin{eqnarray*}     
  \bE_{\sip} (\chi_n)
              & \geq &  (1+o(1)) \, \frac{1}{2} \,  
              \frac{\bE_0 ( \chi_n \, e^{- \sip \, \philone} )}
                           {\bE_0 ( e^{- \sip \, \philone} )} \, 
                      \\*[0.1cm]  
              & \geq &  (1+o(1)) \, m_1(\sip) \, n 
\end{eqnarray*} 
with $ m_1(\sip) = ( \sip b_1)^{1/2} \gamma_1/2 >0,$ possibly depending
on $\sip.$
This ends the proof.
\end{proof} 

Since $\bE Y^2 \geq (\bE Y)^2$ for any random variable $Y,$
we have the following

\begin{corollary}
  \label{vardimone}
Let $\sip > 0$ and $d=1.$ Then there is a constant
$m_2(\sip) >0  $ that may depend on $\sip$
such that as $n \rightarrow \infty,$   
\begin{eqnarray*}
 \bE_{\sip} (\chi_n^2)  & = & (1 + o(1)) \, M_2(n) \, n 
\end{eqnarray*} 
for $ m_2(\sip) \leq  M_2(n) \leq 1.$ 
\end{corollary}


\end{document}